\documentclass{article}

\def\boldpi{{\boldsymbol{\Pi}}}
\def\tct{{n\choose 2}}

\def\OGamma{{\overline{\Gamma}}}
\def\viborita{{\rightsquigarrow}}
\def\hal{\hfill{\vrule height 5 pt width 0.05 in  depth 0.8 pt}\vglue 0.4 cm}
\def\vib{{\rightsquigarrow}}
\def\ogamma{{\overline{\gamma}}}
\def\boldge0{{\boldsymbol{\ge \delta}}} 
\def\OGk{{\overline{G_k}}}
\def\boldle0{{\boldsymbol{\le \delta}}}
\def\boldlt0{{\boldsymbol{< \delta}}}
\def\boldgt0{{\boldsymbol{> \delta}}}
\def\ignore#1{{ }}

\usepackage{ifthen}
\usepackage[paperwidth=7.2in, paperheight=11.5in, width=6.4in,
  height=8.9in, lmargin=0.8in, rmargin=0.6in]{geometry}

\usepackage{amsfonts}
\usepackage{amsmath}
\usepackage{amssymb}

\newtheorem{theorem}{Theorem}[section]

\newtheorem{remark}[theorem]{Remark}
\newtheorem{lemma}[theorem]{Lemma}

\newtheorem{definition}[theorem]{Definition}

\title{Balanced lines in two--coloured point sets}

\author{David Orden\thanks{Departamento de Matem\'aticas, Universidad
    de Alcal\'a, Spain} \and Pedro Ramos\thanks{Departamento de
    Matem\'aticas, Universidad de Alcal\'a, Spain} \and Gelasio
    Salazar\thanks{Instituto de F\'\i sica, Universidad Aut\'onoma de
    San Luis Potos\'{\i}, Mexico}}

\begin{document}

\maketitle

\begin{abstract}
{Let $B$ and $R$ be point sets (of {\em blue} and {\em red} points,
respectively) in the plane, such that $P:=B\cup R$ is in general
position, and $|P|$ is even. A line $\ell$ is {\em balanced} 
if it spans one blue and one red point, and on each open halfplane of
$\ell$, the number of blue points minus the number of red points is
the same.  We prove that $P$ has at least $\min
\{|B|,|R|\} $ balanced lines. This refines a result by Pach and
Pinchasi, who proved this for the case $|B|=|R|$.}
\end{abstract}

\ignore{
\begin{abstract}
{Pach and Pinchasi proved that if $P$ is a point set in general
position in the plane, consisting of $n$ blue points and $n$ red
points, then there are at least $n$ balanced lines, that is, lines
that span one blue and one red point and leave the same number of blue
and red points in each of its open halfplanes. We present a shorter,
simpler proof of this result.}
\end{abstract}
}

\section{Introduction}\label{introduction}

Throughout this discussion, $B$ (respectively, $R$) is a set of 
{\em blue} (respectively, {\em red}) points in the plane, such that
$P:=B\cup R$ is in general position. Let $b:=|B|$ and $r:=|R|$. Without any loss of generality, $b \ge r$.  We assume that $r+b$ is even; equivalently, $r=b + 2\delta$ for some integer $\delta \ge 0$.

A line $\ell$ is {\em balanced} if it spans one blue and one red point, and on each open halfspace of $\ell$, the number of blue points minus the number of red points equals $\delta$. 

It is convenient to assign to each blue point a
{\em weight} of $+1$, and to each red point a weight of $-1$. Thus
a line is balanced if
it spans a point of weight $+1$ and a point of weight $-1$, and the
sum of the weights of the points lying on each of its halfplanes is $\delta$.

Our aim is to find a lower bound for the number of balanced lines.  Pach and Pinchasi~\cite{pachpinchasi} investigated the case $b=r$, and showed that in this case there are at least $b$ balanced lines, and that this bound can be attained (for instance, if $B$ and $R$ are separated by a line).

Our main result is the following.

\begin{theorem}\label{maintheorem}
Let $P$ be a point set in the plane in general position, consisting of
$b$ blue and $r$ red points, with $b \ge r$, and $b+r$ even. Then $P$ has at least $r$ balanced lines.
\end{theorem}

Theorem~\ref{maintheorem} has been indirectly established by Sharir
and Welzl~\cite{sharirwelzl}, who proved that it is equivalent to a
special case of the General Lower Bound Theorem. Our proof is purely
combinatorial.

We emphasize that we have already written a paper with a proof of
Theorem~\ref{maintheorem} using only geometrical
arguments~\cite{ors}. Although the proof in the present paper is
strictly more general that the proof in~\cite{ors} (our proof here is
based on allowable sequences, and so it applies in the more general
setting of generalized configurations of points), the proof
in~\cite{ors} is simpler and more attractive. It would be pointless to
(pretend to) publish both papers, so we have chosen to submit the more
lively, streamlined~\cite{ors} for publication.

\section{Allowable sequences: an overview}\label{section2}

We shall make extensive use of the language and techniques of allowable
sequences, introduced by Goodman and Pollack in~\cite{goodmanpollack}. We now briefly review their definition and basic properties.

Throughout this section, $P$ is as defined in
Section~\ref{introduction}. We define $n:=b + r$.

For simplicity, 
we further assume that no
two lines spanned by two pairs of points of $P$ are
parallel. This is not restrictive at all, since points of $P$ may be
slightly perturbed to satisfy this condition, without altering the
conclusion: a balanced line of the perturbed set corresponds to a
balanced line of the original set.

Let $C$ be a circle that contains all of $P$ in its interior, and let
$\ell$ be a directed line tangent to $C$. We let $\pi^0$ denote the
permutation of $P$ defined by the orthogonal projections of the points of $P$
along $\ell$, read in the direction of $\ell$ (by suitably choosing
$\ell$ we may assume that no two points in $P$ have the same
orthogonal projection over $\ell$).  Now by continuously clockwise
rotating $\ell$, while keeping it tangent to $C$, we eventually find
two points in $P$, say $p$ and $q$, having the same orthogonal
projection. Immediately after that we find that the permutation
defined by the orthogonal projections of the points in $P$ has
changed: $p$ and $q$ have been transposed. We call this new
permutation $\pi^1$. By continuing in this fashion, we define an
infinite sequence $\pi^0, \pi^1, \pi^2, \ldots$ of permutations. By
doing this also in the counterclockwise direction (starting at the
same original line), we also obtain an infinite sequence $\pi^0,
\pi^{-1}, \pi^{-2},\ldots$ of permutations. We combine these into a
doubly infinite sequence $\boldpi:=\ldots, \pi^{-1}, \pi^0,
\pi^1,\ldots$, the {\em allowable sequence} associated to $P$.

A basic property of $\boldpi$ is that two consecutive permutations
$\pi^{t-1}$ and
$\pi^{t}$ differ by a transposition $\tau_{t}$ of two consecutive
points. 
Moreover, for each integer $t$, $\pi^{t+\tct}$ is
the reverse permutation of $\pi^{t}$. Thus $\boldpi$ is periodic
with period $2\tct$, that is,$\pi^{t + 2\tct} = \pi^t$ for
every integer $t$. 

\ignore{
During a rotation of $\ell$ by an angle of $180$ degrees, for each
pair of points of $P$ their
orthogonal projections  along $\ell$ coincide exactly once, say at
angle $\theta$, and as we rotate $\ell$ both counterclockwise and
anticounterclockwise, their orthogonal projections coincide (only) at
all angles $\theta + m\pi$, for $m\in{\mathbb Z}$. It 
}

A transposition $\tau_t: pq\to qp$ is {\em balanced} if (i) one of $p$
and $q$ is blue, and the other is red; and (ii) the sum of the weights
of the points to the left of the pair $pq$ at $\pi^{t-1}$ (and
consequently to the left of the pair $qp$ at $\pi^t$) is $\delta$.  Since the total weight of all points of $P$ is $2\delta$, this implies that the sum of the weights of the points to the
right of the pair $pq$ in $\pi^{t-1}$ (and consequently to the right
of the pair $qp$ at $\pi^t$) is also
$\delta$.

The interest in balanced transpositions is the following key
observation.

\begin{remark}\label{ther}
There is a one--to--one correspondence between balanced
lines of $P$ and the set of balanced transpositions in
its allowable sequence.
\end{remark}

\begin{theorem}[Implies Theorem~\ref{maintheorem}, by
    Remark~\ref{ther}]\label{maintheorems}
If $\boldpi$ is an allowable sequence of $b+r$ points, $b$ blue and
$r$ red, with $b\ge r$ and $b+r$ even, then it has at least $r$ 
distinct balanced transpositions.
\end{theorem}

\section{Proof of Theorem~\ref{maintheorems}}

The proof of Theorem~\ref{maintheorems} relies on the following easy observation: if for some integer $t$ there is a blue point $p$ such that the sum of the weights of the points to the left of $p$ at $\pi^t$ is $\delta$, and the sum of the weights of the points to the left of $p$ at $\pi^{t+1}$ is $\delta-1$, then either (i) $\tau_{t+1}$ moves $p$ to the right, a red point $q$ to the left, and $\tau_{t+1}:pq\to qp$ is balanced; or (ii) $\tau_{t+1}$ moves $p$ to the left, and some other blue point $q$ to the right, in which case obviously $\tau_{t+1}$ is not balanced.

One way around possibility (ii) is not to follow only $p$ from one permutation to the next one, but instead to follow the $k$--th (for some integer $k \ge 1$) element of some blue set $F$, from left to right, at every permutation. Indeed, if it so happens that $p$ is the $k$--th element of $F$ from left to right in $\pi^t$, and $q$ is also in $F$, then if $\tau_{t+1}$ transposed $p$ and $q$ as in (ii), then at $\pi^{t+1}$ $q$ would be the $k$--th element of $F$ from left to right in $\pi^{t+1}$ --- and thus, if we moved our attention from $p$ at $\pi^t$ to $q$ in $\pi^{t+1}$, the sum of the weights to the left of the element in hand would not change: it would be $\delta$ at both times. In conclusion, a ``change of sum of weights" from $\delta$ to $\delta-1$ would always detect a balanced transposition. 

Even though the approach we have just described does not always work (namely, if $q$ is not in $F$), this gives the flavor of the strategy to get our required balanced transpositions.

\subsection{Curves, and the tool to detect balanced transpositions}

A {\em curve} is a mapping $\gamma:{\mathbb Z}\to P$
that assigns to each integer $t$ a point of permutation $\pi^t$, 
and satisfies  $\gamma(t+2\tct) = \gamma(t)$ for every $t\in {\mathbb
  Z}$. 
  
  The term ``curve" is chosen to suggest that, if we visualize $\boldpi$ as a matrix with infinitely many rows and $b+r$ columns (each row being a permutation), then highlighting (i.e., choosing) an element of each permutation draws a curve inside the matrix, especially if from one permutation to the next one the positions of the chosen elements are close to each other. This suggests some concept of continuity, which we will formalize in two variants, namely weak and strong continuity.  

Let 
$Q$ be a nonempty subset of $P$, and let $k$ be an integer, $1 \le k
\le |Q|$. We let $Q_k$ denote the curve that assigns to each integer $t$ the
$k$--th (from left to right) point of $Q$ in $\pi^t$.

For any integer $t$, the {\em weight of $Q_k$ at
time $t$}
is the sum of the weights of the points to the left of $Q_k(t)$ in
$\pi^t$. It is readily checked that $Q_k(t+1)$ and $Q_k(t)$ either coincide or 
are neighbors in $\pi^{t+1}$ (also in $\pi^t$),
and so the weight of $Q_k$ changes by $-1,0$, or $+1$ as we
move from $t$ to $t+1$. We refer to this property as the {\em strong
			continuity} of $Q_k$.

Our (only) tool to find balanced transpositions is the following
statement.%

\begin{lemma}\label{l21a}
Let $F$ be a blue set. Suppose that the weight of $F_k$ is $\delta$ at 
time $t$, and $\delta-1$ at time $t+1$. Then $\tau_{t+1}$ swaps some point $f$ of $F$. Moreover, either (i) $\tau_{t+1}$ moves $f$ to the right and a red point to the left, in which case $\tau_{t+1}$ is 
balanced, and in $\pi^t$ there are exactly $k-1$ points of $F$ to the left of
the transposed points; or (ii) $\tau_{t+1}$ moves $f$ to the left, and moves to the right a point of $B\setminus F$.
\end{lemma}

\noindent{\em Proof.}
It all follows 
from the definition of $F_k$. First, $\tau_{t+1}$ must involve some
point $f$ of $F$, as otherwise the weight of $F_k$ would not change.
Second,  $\tau_{t+1}$
cannot swap two points of $F$, because when this happens the weight
of $F_k$ does not change. Thus $\tau_{t+1}$ swaps $f$ with a point
$x$ not in $F$. 
If $\tau_{t+1}$ moves $f$
to the right, then $x$ cannot be blue, as then the
weight of $F_k$ would increase; thus $x$ is red,
and so $\tau_{t+1}=fx \to xf$ is clearly balanced, and the very
definition of $F_k$ implies that in $\pi^t$ there are exactly $k-1$
points of $F$ to the left of the transposed points. Finally, if $\tau_{t+1}$
moves $f$ to the left, then $x$ cannot be red, for then the weight of
$F_k$ would increase; thus in this case $x$ is in $B\setminus F$. 
\hal

We could say (with good reason) that the proof of the following statement is totally analogous to the proof of Lemma~\ref{l21a}. But a more elegant argument is to simply note that Lemma~\ref{l21b} follows by applying Lemma~\ref{l21a} to the allowable sequence $-\boldpi$ that is the reverse of $\boldpi$, namely $-\boldpi:=(\ldots,\pi^{2},\pi^1,\pi^0,\pi^{-1}, \pi^{-2},\ldots)$. 

\begin{lemma}\label{l21b}
Let $F$ be a blue set. Suppose that the weight of $F_k$ is $\delta-1$ at 
time $t$, and $\delta$ at time $t+1$. Then $\tau_{t+1}$ swaps some point of $F$. Moreover, either (i) $\tau_{t+1}$ moves $f$ to the left and a red point to the right, in which case $\tau_{t+1}$ is
balanced, and in $\pi^{t}$ there are exactly $k-1$ points of $F$ to the left of
the transposed points; or (ii) $\tau_{t+1}$ moves $f$ to the right, and moves to the left a point of $B\setminus F$.
\end{lemma}

A curve {\em is} ${\boldsymbol{\ge
\delta}}$ if its weight is always $\ge \delta$. Similar definitions hold for
$\boldle0$, $\boldlt0$, and $\boldgt0$ curves.

A blue curve (that is, a curve whose image consists only of blue
points) is {\em $\delta$--changing} if it is neither $\boldge0$ nor
$\boldlt0$. A red curve is {\em $\delta$--changing} if it is neither $\boldle0$ nor
$\boldgt0$ (note the assymmetry in the definitions). If a 
monochromatic curve is not
$\delta$--changing, then it is {\em $\delta$--preserving}.

We recall that for every integer $t$,
$\pi^{t+\tct}$ is the reverse permutation of $\pi^t$. This motivates
defining, for each curve $\gamma$, its {\em mirror} curve
$\overline{\gamma}$ by the condition $\overline{\gamma}(t+\tct) =
\gamma(t)$. Thus, if $\gamma$ is blue and its weight at
time $t$ is $r$, then the weight of $\overline{\gamma}$ at time
$t+\tct$ is $2\delta-1-r$. In particular (for blue $\gamma$), whenever the
weight of $\gamma$ is $\ge \delta$, the weight of $\overline{\gamma}$ is
$<\delta$, and viceversa.  

\subsection{Proof of Theorem~\ref{maintheorems} for Case 1: \\  
for every $k$, $\delta + 1 \le k < \lfloor{b/2}\rfloor$, $B_k$ is $\delta$--changing}

We are ready to show Theorem~\ref{maintheorems} for the case in which, for every $k$, $\delta + 1 \le k \le \lfloor{b/2}\rfloor$, $B_k$ is $\delta$--changing, an assumption we make throughout this subsection.

First we claim the following: (A) {\sl For each $k$, $\delta+1\le k \le \lfloor{b/2}\rfloor$, 
			     there are at
least two distinct balanced transpositions that have exactly $k-1$
blue points to the left of the transposed points.}
To see this, first note that since the weight of $B_k$ takes both
$\ge \delta$ and $< \delta$ values, the strong continuity of $B_k$
implies that its weight changes from $\delta$ to $\delta-1$ at least once, and
from $\delta-1$ to $\delta$ also at least once. From Lemmas~\ref{l21a}
and~\ref{l21b}, each such weight change detects a balanced
transposition with the required property (this follows simply since $B\setminus
B=\emptyset$). This proves (A).

Now we claim: (B) {\sl Suppose that $b$ is odd, and let $k_0:=\lceil{b/2}\rceil=
(b+1)/2$. Then there is at least one balanced transposition that has
exactly $k_0-1$ blue points to the left of the transposed points.}
For let $t$ be any integer. If the weight of $B_{k_0}$
at time $t$ is $\ge \delta$ (respectively, $<\delta$), then it is $<\delta$
(respectively, $\ge \delta$) at time $t + \tct$: this follows since (i) the
points to the right of $B_{k_0}$ at time $t$ are precisely the points to 
its left at time $t+\tct$; and (ii) the total weight of all points
different from $B_{k_0}(t)$ is $2\delta-1$. As above, it follows from Lemmas~\ref{l21a} or~\ref{l21b} that either of these weight
changes detects a balanced transposition with the required
property. This proves (B).

We now complete the proof of Theorem~\ref{maintheorems} in case each
$B_k$ is $\delta$--changing. From (A), for each $k, \delta+1\le k\le
\lfloor{b/2}\rfloor$, there are two distinct balanced transpositions
that leave exactly $k-1$ blue points to the left. This provides 
$2(\lfloor{b/2}\rfloor - \delta)$ distinct balanced transpositions. This number
equals $b-2\delta=r$ if $b$ is even, and $b-2\delta-1=r-1$ if $n$ is odd. In the latter case,
(B) yields an
additional balanced transposition. \hal

\subsection{Proof of Theorem~\ref{maintheorems} for Case 2: \\ for some $k$, $\delta+1\le k \le \lfloor{b/2}\rfloor$,  $B_k$ 
 is $\delta$--preserving}

We first outline the strategy to prove Theorem~\ref{maintheorems} if $B_k$ is $\delta$--preserving for some $k$, $\delta+1\le k \lfloor{b/2}\rfloor$. This informal account will help the reader to grasp the flavor of the proof, and will motivate to introduce the ideas and concepts that will be formally defined later.

\vglue 0.4 cm
\noindent{\em An informal account of the strategy}
\vglue 0.4 cm

Suppose that there is a strongly continuous blue curve $\Gamma$  that
is $\boldge0$. Suppose that the curve $\OGamma$ (we recall, the mirror
of $\Gamma$) satisfies that $\Gamma(t)$ is to the left of $\OGamma(t)$
for every $t$. Note that $\OGamma$ is $\boldlt0$.  
We shall call such a $\Gamma$ a border.

We focus on the halfperiod $\pi^0,\pi^1,\ldots,\pi^{\tct}$. Let $F$ be
the set of blue points to the left of $\Gamma$ at $\pi^0$, including
$\Gamma(0)$; let $G$ consist of those blue points between $\Gamma$ and
$\OGamma$ in $\pi^0$; and let $H$ be the set of those blue points to
the right of $\OGamma$ at $\pi^0$, including $\OGamma(0)$.  Thus $B$
is the disjoint union of $F, G$, and $H$.

One main tool to obtain balanced transpositions is the observation
that each time the weight of a curve $F_j$, $1\le j \le |F|$, changes
its weight from $\delta$ to $\delta-1$ (we say this is a $\delta\vib
\delta-1$ change in $F_j$), we detect a balanced transposition: this
follows at once from Lemma~\ref{l21a}, since no blue--blue
transposition in this time interval moves an element of $F$ to the
left, simply because at $t=0$ the first $|F|$ blue elements from left
to right are precisely the elements of $F$.  Similarly, each
$\delta-1\vib\delta$ change in $H_i$, $1 \le i \le |H|$, detects a
balanced transposition.

Now each $F_j$ (respectively, $H_i$) undergoes a $\delta\vib \delta-1$
(respectively, $\delta-1\vib\delta$) change at least once. Indeed,
since $\pi^\tct$ is the reverse permutation of $\pi^0$, and $\Gamma(0)
= \OGamma(\tct)$, it follows that in $\pi^\tct$ each point of $F$ is
either on $\OGamma$ or to its right. Thus each curve $F_j$, $1 \le k
\le |F|$, has to ``cross" both $\Gamma$ (which is $\boldge0$) and
$\OGamma$ (which is $\boldlt0$), and consequently the weight of $F_j$
must change from $\delta$ to $\delta-1$ at least once in this
interval.  An identical argument works for $H_i$: it must change from
$\delta-1$ to $\delta$ at least once in this interval. Thus this
already guarantees at least $|F|+|H|$ distinct balanced
transpositions.
 
 Now we turn our attention to the curves $G_k$, $1 \le k \le |G|/2$
 (suppose for simplicity that $|G|$ is even). Since we are working in
 the time interval between $0$ and $\tct$, it follows that in order to
 recover the full information of a $G_k$ (something we will need), we
 must include $\OGk$ in this time interval. Thus we look at both $G_k$
 and $\OGk$.
  
 By choosing $\Gamma$ appropriately, we will guarantee that for
 $k=1,\ldots,|G|/2$, $G_k\cup \OGk$ undergoes the changes
 $\delta\vib\delta-1$ {\em and} $\delta\vib\delta-1$ in between
 $\Gamma$ and $\OGamma$. Again using Lemma~\ref{l21a}, each
 $\delta\vib\delta-1$ change in $G_k$ or $\OGk$ either detects a
 balanced transposition or induces a $\delta-1\vib\delta$ change in
 some $F_j$ (this crucial technique comes
 from~\cite{pachpinchasi}). Now although one such change in a $F_j$
 does not correspond to a balanced transposition (recall from the
 previous paragraphs that the ones that do are $\delta\vib\delta-1$),
 every time time the weight of $F_j$ changes from $\delta-1$ to
 $\delta$, afterwards it must go down from $\delta$ to $\delta-1$. A
 similar reasoning applies to the curves $H_i$.
 
In conclusion, each $\delta\vib\delta-1$ or $\delta-1\vib\delta$
change in a $G_k$ or $\OGk$ either detects a balanced transposition of
an element of $G$, or indirectly detects a balanced transposition
involving an element $F$ or $H$. Since each $G_k\cup \OGk$, for
$k=1,\ldots,{|G|/2}$ has at least two such changes, this ends up
detecting $2|G|/2= |G|$ balanced transpositions, which added to the
$|F| + |H|$ we had already found, yields $|B|=b$ balanced
transpositions, since $B$ is the disjoint union of $F, G$, and $H$.

Had we started with an analogous red curve $\Gamma$, we would have guaranteed only $r$ balanced transpositions --- which is as good, since this is the number asked for in Theorem~\ref{maintheorems}.

Admittedly, several loose ends remain, among which the most important are: (i) How to guarantee that a $\Gamma$ with all these properties (an appropriate ``border") exists? As it happens, we cannot guarantee the existence of a strongly continuous such $\Gamma$, but of a ``weakly continuous" one; but then, as we will prove, such a $\Gamma$ suffices. (ii) How do we know that each $G_k\cup \OGk$ must have both changes $\delta\vib\delta-1$ and $\delta-1\vib\delta$  in between $\Gamma$ and $\OGamma$? 

With these issues in mind, we now move on to the formal, complete proof of the remaining case of Theorem~\ref{maintheorems}.

\vglue 0.4 cm
\noindent{\em Borders, the partial order $\preceq$ for curves, and $a\vib b$ changes}
\vglue 0.4 cm

If $\alpha, \gamma$ are curves such that for every integer $t$,
$\alpha(t)$ is either equal to or lies to the left of $\gamma(t)$ in
$\pi^t$, then we write $\alpha \preceq \gamma$. If $\alpha(t)$ is
always strictly to the left of $\gamma(t)$, then we write $\alpha
\prec \gamma$. Note that $\preceq$ defines a partial order on the set
of all curves.

\begin{definition}[Borders]\label{borders}
A $\boldge0$ blue curve (respectively, $\boldle0$ red curve) $\gamma$
is a {\em border} if for every integer $t$: (I) either $\gamma(t+1) =
\gamma(t)$, or there are no blue (respectively, red) points between
$\gamma(t)$ and $\gamma(t+1)$ in $\pi^{t+1}$; and (II) $\gamma$ and its mirror curve $\ogamma$ satisfy $\gamma \prec
{\ogamma}$.
\end{definition}

Borders always exist: by assumption, there is a $\delta$--preserving
$B_k$, with $\delta+1\le k \le \lfloor{b/2}\rfloor$. If $B_k$ is $\boldge0$, then it is a border. Otherwise it is
$\boldlt0$, 
in which case it is readily checked that if we define $\rho$
such that for every integer $t$, $\rho(t)$ is the red point in $\pi^t$ to the left of $B_k(t)$ that is
closest to $B_k(t)$ among all red points, then $\rho$ is a
border (that such a red point to the left of $B_k(t)$ exists follows since $k \ge \delta+1$). 

Let $\Gamma$ be a $\preceq$--maximal curve in the set of all borders.

\begin{remark}
For the remainder of the proof, we assume that 
$\Gamma$ is blue. With the assumption of $\Gamma$ being blue, we shall obtain not only $r$, but $b$ balanced transpositions. It is a totally straightforward exercise to adapt every argument to the case in which $\Gamma$ is red, in which case  $r$ balanced transpositions are obtained.    We finally note that from the comment after Definition~\ref{borders} it follows that $B_{\delta+1}\preceq \Gamma$.
\end{remark}

An important property of $\Gamma$ is given in the following statement.

\begin{lemma}\label{unomas}
Let $Q$ be a blue set, and $k$ an integer, $1 \le k \le |Q|$. Suppose
that $Q_k$ is $\boldlt0$.
Then $\Gamma \not\prec Q_k$.
\end{lemma}

\noindent{\em Proof. }
Suppose $\Gamma \prec Q_k$.
Define $\rho$ such
that $\rho(t)$ is the red point in $\pi^t$ to the left of $Q_k(t)$ that is
closest to $Q_k(t)$ among all red points. Such a red point to the left of $Q_k(t)$ always exists, since $Q_k$ is $\boldlt0$, and, moreover, it lies to the right of $\Gamma(t)$, since $\Gamma$ is $\boldge0$. Since $\Gamma$ is blue and $\rho$ is red, then $\Gamma\prec\rho$. It is easy to check that
$\rho$ is a border. Since $\Gamma \prec \rho$, this contradicts the
choice of $\Gamma$.
\hal

A curve $\alpha$ {\em has the change} $a\viborita b$ {\em to the right
  of} curve $\gamma$ if for some time $t$, the weight of $\alpha$ is $a$ at time
  $t$, $b$ at time $t+1$, and in $\pi^t$, $\alpha(t)$ is either
  $\gamma(t)$ or at its right, and in $\pi^{t+1}$, $\alpha(t+1)$ is to the right of
  $\gamma(t+1)$.
Changes {\em to the left} are similarly defined. If $\alpha$ and
  $\beta$ are curves such that $\alpha \prec \beta$, and a change (of
  some curve $\gamma$)
  occurs to the right of $\alpha$ and to the left of $\beta$, then it
  occurs {\em between $\alpha$ and $\beta$}.

If $\alpha$ and $\gamma$ are curves, $\alpha$ {\em crosses $\gamma$
  from the right at time $t$} if $\alpha(t-1)$ either equals or is to the right of
  $\gamma(t-1)$ in $\pi^{t-1}$, and $\alpha(t)$ is to the left of
  $\gamma(t)$ in $\pi^t$. We define the event that $\alpha$ {\em
  crosses $\gamma$ from the left} analogously.
   
A $0\vib -1$ or a $-1\vib 0$ weight change of a blue curve is {\em confined}
if it occurs between $\Gamma$ and $\OGamma$.

\vglue 0.4 cm
\noindent{\em Definitions of $F, G$,and $H$}
\vglue 0.4 cm

We let $F$ denote the set that consists of $\Gamma(0)$ plus 
of all blue points that in $\pi^0$ lie to
the left of $\Gamma(0)$; let $G$ denote the set of all
blue points that in $\pi^0$ lie between $\Gamma(0)$ and $\OGamma(0)$
(including neither of them); and let $H$ denote the set that consists
of $\OGamma(0)$ plus all blue
points that in $\pi^0$ lie to the right of $\OGamma(0)$.
Thus $B$ is the
disjoint union of $F, G$, and $H$.

\vglue 0.4 cm
\noindent{\em Finding confined $\delta\vib\delta-1$ and $\delta-1\vib\delta$ changes in the curves $G_k$ or $\OGk$}
\vglue 0.4 cm

\begin{lemma}\label{gkschangesign}
{\sl Let $k$ be an integer, $1\le k\le\lfloor{|G|/2}\rfloor$. Then $G_k$
has at least one change $0\viborita -1$ and at least one change $-1\viborita 0$ between $\Gamma$
and $\OGamma$. If $n$ is odd, and $k_0:=(|G|+1)/2$, then $G_{k_0}$ has at
least one of the changes $0\vib -1$ or $-1\vib 0$ between $\Gamma$ and $\OGamma$.}
\end{lemma}

\noindent{\em Proof. }
Suppose   $1\le k\le
\lfloor{|G|/2}\rfloor$. 
Suppose that $G_k$ satisfies $\Gamma \prec
G_k\prec \OGamma$ (that is, $G_k$ crosses neither $\Gamma$ nor
$\OGamma$).  We first observe that then $G_k$ cannot be
$\delta$--preserving: for $G_k$ cannot be $\boldlt0$, by
Lemma~\ref{unomas}, and it cannot be $\boldge0$ since it would then be
a border (observe that $\Gamma \prec G_k \prec  \OGamma$), contradicting the choice of $\Gamma$. Thus $G_k$ is
$\delta$--changing. Since $G_k$ (as all curves) has period $2\tct$, it
follows that both changes $-1\vib 0$ and $0\vib -1$ must occur. Thus
we are done if $\Gamma\prec G_k \prec \OGamma$.

Now suppose that $\Gamma\prec G_k \prec \OGamma$ does not hold. Then
$G_k$ must cross $\Gamma$, since in $\pi^0$, $G_k(0)$ lies in between
(and is distinct from both) $\Gamma(0)$ and $\OGamma(0)$.
The periodicity of $G_k$ and  $\Gamma$
implies
it must cross $\Gamma$ both from the right and from the
left.  Now $G_k(0)$ is to the
right of $\Gamma(0)$ at time $0$. Thus let $t_R$ (respectively, $t_L$)
be the smallest
positive (respectively, negative with smallest absolute value) integer
such that $G_k$ crosses $\Gamma$ from the right (respectively, from the
left) at time $t_R$ (respectively, time $t_L$). 
Since $\Gamma \boldge0$, it follows that $G_k$ has weight $< \delta$
at both times $t_L$ and $t_R-1$. Now it cannot happen that $G_k$
has weight $\ge \delta$ at all times between $t_L$ and $t_R-1$, for in this
case we could modify $\Gamma$, making it follow $G_k$ from times $t_L$
until $t_R-1$, while keeping its border properties; %
this would
contradict the choice of $\Gamma$. Thus the weight of $G_k$ is 
$< \delta$ at some time between $t_L$ and $t_R-1$. This implies that
there exist a 
change $0\vib -1$ and a change $-1 \to 0$ to the right of
$\Gamma$. If between $t_L$ and $t_R-1$ $G_k$ does not cross $\OGamma$,
then both changes clearly occur to its left, and so we are done. If in
this time interval $G_k$ does cross $\OGamma$, then, since $\OGamma$
is $\boldlt0$, the first change $0\vib -1$ occurs before the first
such crossing, that is, to the left of $\OGamma$. Thus in either case,
this $0\vib -1$ occurs to the left of $\OGamma$. An analogous
reasoning shows that the first $-1\vib 0$ change before $t_R-1$ also
occurs to the left of $\OGamma$, regardless of whether or not 
$G_k$ crosses $\OGamma$
between $t_L$ and $t_R-1$. 

Finally, if $n$ is odd and $k_0=(|G|+1)/2$, then if the weight of
$G_{k_0}$ at time $0$ is $\ge 0$, then it is $<0$ at time $\tct$, and
viceversa: this follows since $G_k(0)$ is blue and $G_{k_0}(\tct) =
G_{k_0}(0)$ (this last equality holds only for $k_0$). Thus $G_{k_0}$
is $\delta$--changing. Then either $0\vib -1$ or $-1\vib 0$ occur for some
$t$, $0 < t \le \tct$.

An analysis similar to the one in the previous paragraph shows that,
moreover, one such change has to occur between
$\Gamma$ and $\OGamma$.
\hal\medskip

{\em For the rest of the proof we focus exclusively on the time interval
$[0,\tct]$: we work only with objects (such as curves) and events (such as sign changes, or crossings of curves) involving the permutations $\pi^0, \pi^1,\ldots,\pi^{\tct}$, and the transpositions $\tau_1, \tau_2, \ldots, \tau_\tct$.}

\begin{remark}\label{newgame}
It follows from Lemma~\ref{gkschangesign} that 
for
$k=1,\ldots,\lfloor{|G|/2}\rfloor$, the sum of the number of confined 
$0\vib
-1$ changes and confined $-1\vib 0$ changes in $G_k$ or $\OGk$, 
is at least $2$. %
Also, if $n$ is odd and $k_0:=(|G|+1)/2$, then
$G_{k_0}$ has a confined such change.
\end{remark}

\def\ch{{\hbox{\rm ch}}}

\vglue 0.4 cm
\noindent{\em Each confined $\delta\vib\delta-1$ change of a $G_k$ or $\OGk$ either detects a balanced transposition or a $\delta-1\vib\delta$  change in some $F_j$ }
\vglue 0.4 cm

For the following discussion, we assign to every $F_j$, $1 \le j \le
|F|$  a nonnegative integer,
its {\em charge} $\ch(F_j)$. Initially, all
charges are set at $0$.

Consider a confined change $\delta\viborita \delta-1$ of one of the curves $G_k$
or $\OGk$, from some time $t$ to $t+1$. From Lemma~\ref{l21a}, $\tau_{t+1}$ swaps a point $g$ of $G$.  Also by Lemma~\ref{l21a}:

\begin{description}
\item{(a)} If $\tau_{t+1}$ moves $g$ to the right, then $\tau_{t+1}$ is a balanced transposition, and in $\pi^{t}$ there are $k-1$ points of $G$ to the left of $g$, if we are dealing with $G_k$, or $|G|-k$ points of $G$ to the left of $g$, if we are dealing with $\OGk$.  

\item{(b)} 
If $\tau_{t+1}$ moves  $g$ moves to the left, then $\tau_{t+1}$ moves to the right a point $b$ of $B\setminus G$; moreover, in our working time interval, $b$ must belong to $F$, since a swap involving a point of $G$ and a point of $H$ move the point of $G$ to the right.  If $f$ is the $j$--th point
of $F$ from left to right in $\pi^t$, then $\tau_{t+1}$ induces a 
weight change $\delta-1\viborita \delta$ of curve $F_j$. In this case we let this
transposition increase the charge of $F_j$ by $1$ (this is a {\em
charge transaction}).

\end{description}

\vglue 0.4 cm
\noindent{\em Each confined $\delta-1\vib\delta$ change of a $G_k$ or $\OGk$ either detects a balanced transposition or a $\delta-1\vib\delta$ change in some $H_i$}
\vglue 0.4 cm

As we did with each $F_j$, we now assign to every $H_i$, $1 \le i \le |H|$ a nonnegative integer,
its {\em charge} $\ch(H_i)$. Initially, all
charges are set at $0$.

An totally analogous argument to the one above shows that a confined change $\delta-1\vib \delta$ of one
of the curves $G_k$ or $\OGk$ either reveals: {(a')} a balanced
transposition, involving a point of $G$ and leaving $k-1$ points
of $G$ to its left, if we are dealing with $G_k$, and $|G|-k$
points of $G$ to its left, if we are dealing with $\OGk$; or {(b')} a
change $\delta\vib \delta-1$ of a curve $H_i$. In the latter case, we let this
transposition increase the charge of $H_i$ by $1$ (also a
{\em charge transaction}).

\vglue 0.4 cm
\noindent{\em Bounding (by below) the number of $\delta\vib \delta-1$ changes in curves $F_j$ and the number of $\delta-1\vib \delta$ changes in curves $H_i$}
\vglue 0.4 cm

\noindent {\bf Claim A } 
{\sl For each $j\in\{0,\ldots,|F|\}$, the number of times
that curve $F_j$ has the change $\delta\viborita \delta-1$ (in our working time interval) is at least 
$\ch(F_j) + 1$.}

\noindent{\em Proof. }
First note that $\Gamma$ is a
$\boldge0$ blue curve, and so just before $F_j$ first crosses $\Gamma$ 
($F_{|F|}$ may never cross $\Gamma$, but then at $t=0$ its weight is $\ge \delta$) its
weight is $\ge \delta$; afterwards, it undergoes at least $\ch(F_j)$
confined changes of the type $\delta-1\viborita \delta$.
 These together imply at least $\ch(F_j)$ changes of the
type $\delta\viborita \delta-1$ for $F_j$, each of them occurring {\em before} an
identified $\delta-1\viborita \delta$ change. Now {\em after} the last
confined $\delta-1\viborita \delta$ change, we must have an additional $\delta\viborita \delta-1$
change, since whenever
$F_j$ last crosses $\OGamma$, the weight of $F_j$
is $<\delta$. 
\hal\medskip

A totally analogous argument shows the following.

\noindent {\bf Claim B } 
{\sl For each $i\in\{0,\ldots,|H|\}$, the number of times
that curve $H_i$ has the change $\delta-1\viborita \delta$ (in our working time interval) is at least 
$\ch(H_i) + 1$.}

\vglue 0.4 cm
\noindent{\em Concluding the proof}
\vglue 0.4 cm

We observe that Lemma~\ref{l21a} applied to (each) $F_j$ implies that {\em each $\delta\vib \delta-1$ change corresponds to a balanced transposition involving a point of $F$ and leaving $j-1$ points of $F$ to the left}. This follows since conclusion (ii) in Lemma~\ref{l21a} cannot apply here, since in our working interval no point of $f$ moves to the left transposing with a point of $B\setminus F=G\cup H$. 

Similarly, Lemma~\ref{l21b} applied to (each) $H_i$ implies that {\em each $\delta -1 \vib \delta$ change corresponds to a balanced transposition involving a point of $H$ and leaving $i-1$ points of $H$ to the left}.

The previous two paragraphs, combined with Claims A and B, imply that  the number of distinct balanced transpositions involving a point of $F$
or $H$ 
is at least $\sum_{j=1}^{|F|} (\ch(F_j) + 1)$ 
$\sum_{i=1}^{|H|} (\ch(H_i) + 1) = |F| + |H| + 
\sum_{j=1}^{|F|} \ch(F_j) + 
\sum_{i=1}^{|H|} \ch(H_i)$. Therefore:

\noindent{(I)} {\sl the number of distinct balanced
transpositions involving a point of $F\cup H$ is at least $|F| +
|H|$ plus the total number of charge transactions.}

Now for $k=1,\ldots,\lfloor{|G|/2}\rfloor$, there are at least $2$
confined changes $\delta\vib \delta -1$ or $\delta -1\vib \delta$ involving $G_k$ or
$\OGk$. If $|G|$ is odd and $k=(|G|+1)/2$, then there is at least $1$
such confined change. Now we recall (see (a), (b), (a'), and (b')) that each such change 
identifies either a balanced transposition of a point of $G$ or
results in a charge transaction. The balanced transpositions
identified are all distinct, since a transposition identified
following a $G_k$ (or $\OGk$) leaves $k-1$ (or $|G|-k$) points of
$G$ to its left.
Therefore:

\noindent{(II)} {\sl 
The number of distinct balanced
transpositions involving a point of $G$ is at least
$2(\lfloor{|G|/2}\rfloor)$ ($+1$, if $|G|$ is odd) minus the total
number of charge transactions. }

Combining (I) and (II), we conclude
that, regardless of the parity of $|G|$, the total number of distinct
balanced transpositions is at least $|F| + |G|+ |H|$, which equals
$|B|=b$, since $B$ is the disjoint union of $F, G$, and $H$.

\end{document}